\documentclass[10pt]{article}
\usepackage{amssymb,amsmath,amsfonts} 
\usepackage[colorlinks=true]{hyperref}
\usepackage{comment}

\newtheorem{remark}{Remark}

\newcommand{\io}{{\infty}}

\newcommand{\real}{ {\mathbb R}   }
\newcommand{\torus}{ {\mathbb T}   }
\newcommand{\integer}{ {\mathbb Z}   }
\newcommand{\complex}{ {\mathbb C}   }

\newcommand{\cB}{ {\mathcal B}   }
\newcommand{\cE}{ {\mathcal E}   }
\newcommand{\cP}{ {\mathcal P}   }
\newcommand{\cM}{ {\mathcal M}   }

\newcommand\beq[1]{ \begin{equation}\label{#1} }
\newcommand{\eeq}{ \end{equation} }

\newcommand\beqa[1]{ \begin{eqnarray} \label{#1}}
\newcommand{\eeqa}{ \end{eqnarray} }
\newcommand{\beqano}{ \begin{eqnarray*} }
\newcommand{\eeqano}{ \end{eqnarray*} }

\newtheorem{definition}{Definition}[section]
\newcommand\dfn[1]{ \begin{definition}\label{#1} \rm}
\newcommand\edfn{ \end{definition} }

\newcommand\equ[1]{{\rm (\ref{#1})}}
\newcommand{\nl}{{\smallskip\noindent}}
\newcommand{\giu}{{\medskip\noindent}}
\newcommand{\Giu}{{\bigskip\noindent}}

\newcommand{\x}{\xi}

\newcommand{\e}{\varepsilon}

\renewcommand{\b }{\beta }

\newcommand{\ii }{{\rm i} }
\renewcommand{\d }{\delta }

\newcommand{\f }{\varphi}

\newcommand{\m }{\mu }

\renewcommand{\t }{\tau }
\renewcommand{\o }{\omega }

\newcommand{\Z}{\mathbb{Z}}

\title{On the measure of Lagrangian invariant  tori in nearly--integrable mechanical systems \\
\large (Draft)
}

\begin{document}

\author{ 
\footnotesize L. Biasco  \& L. Chierchia
\\ \footnotesize Dipartimento di Matematica e Fisica
\\ \footnotesize Universit\`a degli Studi Roma Tre
\\ \footnotesize Largo San L. Murialdo 1 - 00146 Roma, Italy
\\ {\footnotesize biasco@mat.uniroma3.it, luigi@mat.uniroma3.it}
\\ 
}

\maketitle


\begin{abstract}
\noindent
Consider a real--analytic nearly--integrable mechanical system with potential $f$, namely,  
a Hamiltonian system, having a real-analytic Hamiltonian  
\beq{H}
H(y,x)=\frac12 | y |^2 +\e f(x)\ ,
\eeq
 $y,x$ being $n$--dimensional standard action--angle variables (and $|\cdot|$ the Euclidean norm).
Then, for ``general'' potentials $f$'s and  $\e$ small enough, the Liouville measure of the complementary of invariant  tori is smaller than $\e|\ln \e|^a$ (for a suitable $a>0$).
\end{abstract}

\section{Introduction and results}

The main result of last century in the theory of nearly--integrable Hamiltonian systems is that, under suitable non--degeneracy and regularity assumptions, ``most'' of the regular solutions of the integrable regime, which span Lagrangian tori in phase space, persist under small perturbations. This celebrated  result is due to N.N.~Kolmogorov (1954, \cite{K}) and is the core of the so--called KAM (Kolmogorov--Arnold--Moser) theory (compare \cite{AKN} and references therein). ``Most'' , as it was later clarified in \cite{L}
(in dimension 2) and
\cite{Nei},
\cite{poschel1982}   (in any dimension),  means that  the union of  Lagrangian invariant tori on which the flow is conjugated to a Diophantine linear flow\footnote{A ``Diophantine linear flow'' on the $n$--dimensional standard torus $\torus^n:=\real^n/(2\pi\real)^n$ is the flow $t\to x_0+\o t$ where $\o\in\real^n$ is a Diophantine vector, i.e., satisfies $|\o\cdot k|\ge \kappa/|k|^\t$ for all integer non--null vectors $k\in\integer^n$, for some $\kappa>0$ and $\t\ge n-1$; compare, e.g., \cite{AKN}.}, has relative
measure of order $1-\sqrt{\e}$, if $\e$ is the perturbation parameter. This statement is sometimes rephrased by saying that the ``non--torus (invariant) set''  is of measure less than $\sqrt{\e}$. Indeed, it is easy to see that this result is optimal if one consider only Lagrangian tori which are graphs over the angle variables (hence, are homotopically non--trivial).
In fact, just consider the simple pendulum  $\frac12 y^2 +\e \cos x$ with $y\in \real$ and $x\in\torus^1$: the phase region inside the separatrix (corresponding to oscillations of the pendulum) does not contain any invariant torus (circle) which is a global graph over the angle on $\torus^1$ and this region has measure $4\sqrt{2\e}$.
Of course, in this trivial integrable example, the full phase space is covered by Lagrangian invariant  tori with the exception of the stable equilibrium and the unstable (hyperbolic) equilibrium together with the separatrix (which coincides with the stable and unstable one--dimensional manifold of the hyperbolic equilibrium).
The point is that   invariant tori may have different topologies: in the pendulum case,
the tori corresponding to librations (circles above or below the separatrix)  are homoptopically non trivial, while the tori corresponding to oscillations (``secondary tori''), are contractible circles enclosed by the separatrix.
The tori  obtained by standard applications of  KAM theory in mechanical systems  are homotopically non--trivial and are deformations of the integrable tori  $y=$const. 

\nl
In general, when $\e\neq 0$ ``secondary'' (homotopically different) invariant Lagrangian tori arise near ``resonances'', i.e., regions of phase space where $\omega\cdot k=0$, 
where $\omega$ is the unperturbed frequency, 
``$\cdot$''  the standard inner product and $k$ an integer non vanishing vector (in the pendulum example $\omega = y=0$ and  $k=1$). It is therefore natural to expect that the non--torus region is of measure, in general,  {\sl smaller} than $\sqrt{\e}$. On the other hand   a simple argument (see Remark~(ii) below) suggests that the non--torus region is, in general, of measure larger than $\e$. 
Indeed, Arnold, Kozlov and Neishtadt  conjecture that such region has measure of order $\e$ \cite[Remark 6.18]{AKN}.

\nl
{\sl We can prove that Arnold, Kozlov and Neishtadt's  conjecture is ``essentially''  true in the case of mechanical systems}.

\nl
The word ``essentially'' means that the result holds up to logarithmic corrections in $\e$ and  for a suitable class of real--analytic potentials, which is of full measure  (on a natural probability function space), contains an open dense set and is prevalent\footnote{A Borel set $P$ of a Banach space $X$ is called {\sl prevalent} if there exists a compactly supported probability measure $\nu$ such that $\nu(x+P)=1$ for all $x\in X$; compare, e.g., \cite{HK}.\label{nurzia}}. 

\nl
A precise statement and a sketch of proof will be given below   
(full proofs  will appear elsewhere \cite{BC}).

\noindent
{\bf Functional setting, probability measures and Fourier--projections} 

\nl 
We proceed to describe the ``good'' set of potentials $f$ in \equ{H}, for which the  result holds. 
Roughly speaking,  such  a set  consists of real--analytic functions\footnote{It would be easier to consider larger function spaces of smooth functions. However, the natural (both from the theoretical and applicative point of view) and most challenging setting is definitively that of real--analytic potentials.}, whose Fourier--projections on one--dimensional resonance vectors have a derivative which is a   Morse functions. 

\nl
Let $s>0$ and 
consider the real--analytic functions on $\torus^n$ having zero average and finite norm
$$
|f|_s:=\sup_{k\in \integer^n} |f_k| e^{|k|s}<\infty
$$
where $f_k$ denotes Fourier coefficients and, as usual, $|k|$, for integer vectors, denotes the 1-norm $\sum |k_j|$. Denote by $ \cB_s^n$ the Banach space of such functions. \\
Now, let $\integer^n_\sharp$ denote the set of integer vectors $k\neq 0$ such that the  first non--null  component is positive, 
and denote 
$\ell_\io^n$  the Banach space of complex sequences $z=\{z_k\}_{k\in \integer^n_\sharp  }$ with finite sup--norm 
$|z|_\io:=\sup_{k\in\integer^n_\sharp } |z_k|$.  
The map 
\begin{equation}\label{miserere}
j:f\in  \cB_s^n \to \big\{f_k e^{|k|s}\big\}_{k\in \integer^n_\sharp}\in \ell_\io^n
\end{equation}
is an isomorphism of Banach spaces\footnote{Recall that since the functions in $\cB^n_s$ are {\sl real}--analytic one has the reality condition $f_k=\bar f_{-k}$.},
which allows to identify functions in $ \cB_s^n$ with points in $\ell_\io^n$ and the Borellians of $ \cB_s^n$ with those of $\ell_\io^n$. 
Now, consider the standard normalized Lebesgue--product measure 
on the unit closed ball of $\ell_\io^n$, namely,   the unique probability measure $\mu$ on the Borellians of $\{z\in \ell_\io^n: |z|_\io\le 1\}$ such that,
given   Lebesgue measurable  sets in the unit complex disk $E_k\subseteq D:=\{w\in\complex:\ |w|\le 1\}$
with $E_k\neq D$ only for finitely many $k$, one has
 $$
\mu \Big(\prod_{k\in\Z^n_\sharp} E_k\Big)=
\prod_{\{k\in \integer^n_\sharp:\, E_k\neq D\}} 
{\rm meas}(E_k)\,
$$
where ``meas'' denotes  the  normalized Lebesgue measure on the unit complex disk $D$. 
{\sl The isometry $j$ in \eqref{miserere}
naturally induces a  measure
$\mu_s$ on the
unit ball of $\cB^n_s$}.

\nl
Now,  denote by ${\integer}^n_*$ the set of  vectors $k\in  \integer^n_\sharp$ such that the greater common divisors of their components is 1. Then, any function $f\in  \cB_s^n$
can be uniquely decomposed, in Fourier space, as sum of real--analytic functions of one--variable, which are the projection of $f$ onto the one--dimensional Fourier modulus $\{jk:j\in \integer\}$, as follows\footnote{Beware of notation: $\xi\to F_k(\xi)$ is a periodic function of one variable, whose Fourier coefficients, for $j\in \integer$, are given by $(F_k)_j=f_{jk}$.}:
\beq{dec}
f(x)= \sum_{k\in \integer^n_*} F_k(k\cdot x)\ ,\quad {\rm where}\quad
F_k(\xi):=\sum_{j\in \integer\backslash\{0\}} f_{jk}e^{ij\xi}
\eeq
$f_{jk}$ being the Fourier coefficient of $f$ with Fourier index $jk\in\integer^n$. 
Notice that, since $f\in  \cB_s^n$, the functions $F_k$ belong to $ \cB_{|k|s}^1$. 

\Giu
{\bf The class $\cP_s$ of ``good'' potentials } 

\nl
We first define  an auxiliary function $K_s(\d):=\frac{2}{s}\ln\frac{c}{\d}$ (where $c>1$ depends only on $n$).

\nl 
 {\sl Let $\cP_s$  be the set of functions in the closed unit ball of $ \cB_s^n$ 
satisfying the following three properties. There exists $\d>0$ such that, for $k\in\integer^n_*$,

\giu
{\rm (P1)}
$\displaystyle{  |f_k|\geq \d |k|^{-\frac{n+3}2}\ e^{-|k|s} \,,\ \ 
\forall\  |k|> K_s(\d)\,;}$

\giu
{\rm (P2)}
$\displaystyle{\  
\min_{\xi\in \real} \ \big( |F_k'(\xi)|+|F_k''(\xi)|\big) >0 \,,\ \ 
\forall\  |k|\le  K_s(\d)\,;}$

\giu
{\rm (P3)} $3 F_k''(\bar \xi)  F_k''''(\bar \xi)\neq 5  (F_k'''(\bar \xi))^2$,
for $\bar \xi$ minimum  of  $F_k$ and  $\forall\  |k|\leq  K_s(\d)$.}

\giu
Fix $s>0$ and fix a bounded region $B\subset\real^n$. Then one has the following 

\Giu
{\bf Theorem}\\ {\sl 
{\bf (a)} Let $f\in \cP_s$ and $H:= \frac12 |y|^2 +\e f(x)$.
There exist $\e_0>0$ and $a>0$ such that, for any $0<\e<\e_0$,  the measure of the set of $H$--trajectories in $B\times \torus^n$, which do not lie on an invariant Lagrangian (Diophantine) torus, is bounded by  $\e |\ln \e|^a$.
\\
{\bf (b)} The set $\cP_s$ has full $\m_s$--measure,
contains an open dense set (in the unit ball of  $ \cB_s^n$) and is prevalent. 
}

\Giu
{\bf Comments and remarks} 

\begin{itemize}

\item[(i)] We stress that the set $\cP_s$ is ``large'' in many ways: it is of full measure with respect to a quite natural product probability measure on a weighted Fourier space;
it is generic in the topological sense (Baire), and actually is more than that, since it contains an open dense set; finally, it is prevalent, which is a measure--theoretic notion for subsets of infinite--dimensional spaces that is analogous to ``full Lebesgue measure'' in Euclidean spaces (compare \cite{HK}).

\item[(ii)]
Let us give a simple heuristic argument suggesting that  the measure of the non--torus set is, in general, at least
$O(\e)$. Let $n=2$,  let $p=\sqrt{\e}y$, $q=x$ and divide the Hamiltonian \equ{H} by $\e$: we get a new Hamiltonian, $\tilde H=|p|^2/2 + f(q)$, which is parameter free\footnote{Trajectories of  $H$ and $\tilde H$ are in 1-1 correspondence (but with different times).}. In general, one expects an $O(1)$ torus--free region around in the $(p,q)$ variables: this  corresponds to a torus--free region of measure $O(\e)$ in the original variables. When  $n\ge 2$ the same argument applies in neighborhoods of double resonances and lead to the same conclusion. 

\item[(iii)] As far as we know,   the only other general result  about the measure of secondary tori nearby simple resonances  is discussed by
Medvedev, Neishtadt and  Treschev  in \cite{MNT}. They proved
 that there is a set $\mathcal D_0$ inside the separatrix ``eye'' $\mathcal D$ arising near a simple resonance with meas$\,(\mathcal D_0)\sim\sqrt\e$ ($\sim$ meas$\,(\mathcal D)$), such that  the non--torus set in $\mathcal D_0$  is of measure $O(\e)$.
 This, of course, does not imply that the non--torus set in bounded regions is $O(\e)$, which is what is stated in the Theorem above.

\item[(iv)] We are interested here in global statements,  which apply to any region in phase space. On the other hand, it is well known that the density of invariant  tori is non uniform in phase space and there are regions where such density is of order
$1- e^{-1/\e^c}$ for some $c>0$; compare the Remark at p. \pageref{spuntature}.

\item[(v)]
Even though the proof  of the Theorem  -- briefly sketched below -- exploits the particular form of the Hamiltonian $H$ in \equ{H}, we believe that extending the proof to more general settings  such as to Hamiltonians of the form $H=h(y)+\e f(y,x)$ with $h$ quasi--convex is an interesting but essentially technical matter. \\
A more challenging question would be to determine  the most general class of unperturbed Hamiltonians $h$ for which the result holds. 

\item[(vi)] Two naive questions.  Can one take out the logarithms? 
Is the result true for {\sl all} potentials?

\item[(vii)] We close this section with a technical comment on the definition of $\cP_s$.\\
Condition (P1) insures that, for 
$|k|> K_s(\d),$ the first Fourier coefficient in the expansion of $F_k,$
namely $(F_k)_1=f_k,$ is large with respect to the other ones, since
$|(F_k)_j|=|f_{jk}|\leq e^{-|j||k|s}.$
Roughly speaking, this means that  $F_k$ ``behaves like a cosine''.\\
(P2) says that $F_k'$ is a  ``Morse function''
uniformly for $|k|\leq K_s(\d)$. In particular, 
all critical points of  $F_k$ 
are non degenerate minima or maxima and  the number of critical points of $F_k$ is bounded by a constant times
$\beta^{-1}$, where $\beta$ is the strictly positive  number in the left hand side of (P2). 
\\
(P3) is needed  to prove a non degeneracy condition
for the action-angle variables of one--dimensional systems of the type
$\eta^2/2+F_k(\xi).$

\end{itemize}

\section{Outline of the proof of the Theorem}

\nl 
We begin by discussing {\bf part (a)}, whose proof may be divided in five steps.

\nl
{\sl Step 1. Geometry of resonances}

\noindent
Define
$$ 
\Z^n_\e:=\{\ k\in\Z^n_* \ {\rm s.t.}\ 0<|k|\leq |\ln \e|^2 \ \} 
$$ 
and, for $k\in\Z^n_\e,$
$$
R_k:=\{y\in B: |y\cdot k| \leq  \sqrt{\e}|\ln\e|^{c_n}\}\,,
$$
with $c_n:=2n+6.$

\nl
Decompose the phase space $\cM=B\times\torus^n$ as\footnote{The symbol $\sqcup$ denotes disjoint union.}: 
$B=B_0\  \sqcup B_1\ \sqcup B_2$, where:
 
$$  
B_0  :={B\setminus \bigcup_{k\in\Z^n_\e}} R_k\ , \quad
B_2 :=   \bigcup_{k\neq k'\in \Z^n_\e}   R_k\cap R_{k'}\ , \quad
B_1 := \bigcup_{k\in\Z^n_\e} R_k\setminus B_2\ .
$$
Roughly speaking (and for $|k|\leq |\ln 	\e|^2$), $B_0$ is the non--resonant set,
$B_1$ contains simple resonances and $B_2$ contains double (or higher) resonances.

\noindent
On the non--resonant set $B_0$, by an averaging procedure,
we can remove the perturbation at order $\e^{|\ln \e|}.$
Then, applying a standard KAM theorem we obtain that 
the measure of the non--torus set in $B_0$ is 
$O(\e^{|\ln \e|/2})\ll \e.$ 

\nl
The set $B_2$ can be disregarded since it is of measure $O(\e|\ln \e|^{2 c_n+2})$.  

\nl
The main problem comes from the neighborhood $B_1$ of  simple resonances   (a set whose measure  is not negligible since {meas}$(B_1)\geq\sqrt \e\gg\e$).

\nl
{\bf Remark.} \label{spuntature}
The geometry of resonances here
is different from the geometry of resonances (in the convex case) as discussed, e.g.,  in  \cite{poschel}. In fact, in \cite{poschel}
more resonances are disregarded in the non-resonant set,
 namely, the resonances with  $|k|\leq (1/\e)^a$, $a>0$; moreover the neighborhood of  simple resonances  has width $\e^b,$ $0<b<1/2$, which is larger than $R_k$.
As a consequence, the set of double resonances has measure greater than
$\e^{2b}$, which is a set not negligible for our purposes.
On the other hand in  Nekhoroshev's
theorem one can average out the perturbation up to an exponentially
high order $e^{- {\rm const\, }(1/\e)^a},$ while we get only $\e^{|\ln\e|}$
(see \eqref{vanitatum} below).

\giu
{\sl Step 2. 
Averaging over a simple resonance} 

\nl
Fix $R_k$ with $k\in\Z^n_\e$. After performing

\nl
 
- a linear change of variables putting the resonance 
in $\{y_n=0\},$

- an averaging over the ``fast angles'' $x_1,\ldots,x_{n-1},$

- a $1/\sqrt\e$--blow up of $y$ and an $\e$-time rescaling,  
 
\nl 
the Hamiltonian takes the form: 
\begin{equation}\label{vanitatum}
h(\hat y)+ \frac{1}{2} y_n^2 + \tilde F_k (y,\|k\|^2 x_n) + \e^{|\ln \e|} g_k(y,x)\,,\qquad\quad
\hat y:=(y_1,\dots,y_{n-1})\,,
\end{equation}
where 
\beq{cuordimela}
\tilde F_k(y,\|k\|^2 x_n) = F_k(y_n,\|k\|^2 x_n) +O(|\ln \e|^{-c})\ ,
\eeq
 $F_k$ being defined in \eqref{dec} and satisfying, for some $\d>0$, 
(P1)--(P3), (and $\|k\|^2:=\sum_{j=1}^n k_j^2$). 
Notice that the Hessian of $h$ is non degenerate since $h$ is convex.
Note also that the unperturbed Hamiltonian in \eqref{vanitatum},
 i.e., when $\e=0,$ is $2\pi/\|k\|^2$-periodic in $x_n,$ while the complete Hamiltonian
 is only $2\pi$-periodic in $x_n$, a well known fact due to the presence of a 
 $k$-resonance.

\nl
To simplify the exposition, in view of \equ{cuordimela}, henceforth we replace in \equ{vanitatum} $\tilde F_k$ with $F_k$.
 
\giu
{\sl Step 3. 
Action--angle variables}

\nl
We want  to prove that

\nl
$(*)$\label{magnificat} the measure of the non--torus
set in a $O(|\ln \e|^{c_n})$-neighborhood of $\{y_n=0\}$
for 
the Hamiltonian 
in \eqref{vanitatum}
is $O(\sqrt\e |\ln\e|^c)$,

\nl
which, in view of  the $1/\sqrt\e$--blow up of the action variables (Step 2), 
corresponds, in the original variables, to a  measure of $O(\e |\ln\e|^c).$

\nl
The idea is to use action-angle variables
$(p_n,q_n)\mapsto(y_n,x_n)$ to integrate the one degree--of--freedom
Hamiltonian 
$$
 {E_k(y_n,x_n):=\frac{1}{2} y_n^2 +  F_k (y_n,\|k\|^2 x_n)=
E_k(p_n).}
$$
Thus, in the new variables (completed with $\hat p:=\hat y,$ $\hat q:=\hat x$)
the full Hamiltonian \equ{vanitatum} becomes
\begin{equation}\label{laudateomnesgentes}
h(\hat p)+  {E_k(p_n)}	+\e^{|\ln\e |} \tilde g_k(p,q)\,,
\end{equation}
a form suitable (in principle) for applying a standard KAM theorem. 
However, there are two (quantitative) problems to overcome:
\\ 
{\bf (S)} the transformation  putting $E_k$ in action-angle variables becomes 
singular as separatrices or elliptic equilibria are approached;
 \\
{\bf (K)} Kolmogorov's non--degeneracy condition (namely that the Hessian of 
$h(\hat p)+  E_k(p_n)$ is invertible) needs to be checked.
 
\nl
These two problems will be overcome in the next two steps. 
 
\giu
{\sl Step 4.  On the singularities of action-angle variables  }
 
\nl
This analysis will be based on the following lemma:

\nl
{\sl Assume that in the one--dimensional system $E(\eta,\xi):=\eta^2/2 + F(\xi)$
the derivative of the potential $F(\xi)$ is a Morse function with 
$\b:=\min_{\xi\in \real} \ \big( |F_k'(\xi)|+|F_k''(\xi)|\big)>0$.
Then, there exists $c>0$, depending only on $\b$,  such that
 for any critical energy\footnote{Namely $E_0:=F(\xi_0)$ with $F'(\xi_0)=0$ and, therefore,
the level set ${E=E_0}$
contains equilibria and/or separatrices.} 
$E_0$
the measure of the points $(\eta,\xi)$ for which  $E(\eta,\xi)$ is $\theta$-close to $E_0$
is bounded by $c\,  \theta|\ln \theta|.$}

\nl
We will discuss, here,  only modes $k$ for which   $|k|\leq K_s(\d)$,
since the case $|k|> K_s(\d)$ is easier (recall that assumption (P1) implies  that $F_k$ ``behaves like a cosine''). 

\nl
Now,  we can use the above result thanks to  (P2) with 
\beq{ventaglietto}
\theta=\e^{a |\ln \e |}\ ,\quad 0<a<a_0
\eeq
for a suitable $a_0<1$.
Thus, up to a small set of measure $\theta |\ln \theta|$, we can use ``non--singular'' action--angle variables $(p_n,q_n)$. In particular, the $p_n$--radius of analyticity , the width of the complex $q_n$--strip 
and the $C^2$--norm of $E_k$ are of order, respectively, $\theta$,
$1/|\ln \theta|$ and $1/\theta$.

\giu
{\sl Step 5. On Kolmogorov's condition}
 
\nl
Since $h$ is (strictly) convex, Kolmogorov's non--degeneracy condition amounts to say that
$|E_k''(p_n)|$ is bounded away from zero  and up to a small measure  set of $p_n$.

\nl
We also notice that for  ``high energies'', namely  for energies larger
than then the maximal critical energy, 
$E_k(p_n)$ is strictly convex   and the Kolmogorov's condition is easily satisfied.
Therefore, we need only to discuss the case inside separatrices where $E_k''$ may become negative or  null.

\nl
As above, by (P1), the modes $|k|> K_s(\d)$ 
can be easily handled directly.

\nl
The case $|k|\leq K_s(\d)$ is more difficult:  although it involves only a finite number of $k$'s, the structure of the Hamiltonian $E_k(p_n)$ is rather ``arbitrary'' and it is not at all obvious how to handle it in a direct way. 
To overcome this problem we will
 check the non degeneracy condition in an {\sl  indirect way}, using
the {\sl analyticity} of the function $E_k$.

\nl
$F_k$ has a finite number of critical points, which are  non--degenerate local maxima or minima  (recall the comment (vii) on  (P2)).  The critical energy levels determine a finite number of open   connected components, where one can define analytic action-angle variables. Recall that we are discussing  the {\sl bounded} components and
let us fix one of such components. Let the action variable $p_n$, correspondingly,  be defined on an interval $(a,b)$ on which $E_k(p_n)$ will be a real--analytic strictly increasing function.  
By construction $e_a:=\lim_{p_n\to a^+}E_k(p_n)$ and $e_b:=\lim_{p_n\to b^-}E_k(p_n) $ correspond to critical 
values of $ F_k$.
It is simple to see that $e_b$  corresponds to a (local) maximum value of 
$F_k,$ while $e_a$  corresponds to a (local)  maximum or minimum value of 
$ F_k$.
It is also simple to see that 
$\lim_{p_n\to b^-}|E_k''(p_n)|=\infty$
(and the same holds for $\lim_{p_n\to a^+}|E_k''(p_n)|$  if also $e_a$  corresponds to a maximum).
Moreover, when $e_a$  corresponds to a minimum, 
we can prove that
$\lim_{p_n\to a^+}E_k''(p_n)\neq 0$
using assumption (P3).
In any case we have  
$\lim_{p_n\to a^+}E_k''(p_n)\neq 0,$ $\lim_{p_n\to b^-}E_k''(p_n)\neq 0. $
Then,
thanks to a basic  property of level sets of analytic 
functions\footnote{
{\sl Let $f$ be a non identically zero (real) analytic function in an open interval containing the closed interval $[x_1,x_2]$. 
Then there exists a constant $c\geq 1$ such that for every 
$0<\theta<1/c$ 
$$
{\rm meas}\,  \{ x\in[x_1,x_2] :\ |f(x)|\leq \theta \}  
\leq
c\, \theta^{1/c}\,.
$$}}, 
we conclude that there exists $c>1$ such that, up to a set of $p_n$'s of measure $\theta$ (as in \equ{ventaglietto}), 
the estimate $|E_k''|\geq \theta^c$ holds.
Since the size of the perturbation in \eqref{laudateomnesgentes}
is $O(\e^{|\ln\e|})$ this is enough to apply a KAM theorem (recall the choice of parameters given at the end of Step 4 and choose $a$ small enough in \equ{ventaglietto}). This concludes the sketch of the proof of part {\bf (a)}.

\giu
Let us now outline the proof of  part {\bf (b)}.

\nl
To check that $\m_s(\cP_s)=1$ we shall prove that, for every $\d>0,$ 
the measure of the sets of
potentials $f$ that do not satisfy, respectively, (P1), (P2), (P3) is,
respectively,  $O(\d^2),0,0$; the result will follow letting $\d\to 0.$

\nl
First, by the identification \eqref{miserere}, the
measure of the set of
potentials $f$ that do not satisfy (P1)
with a given $\delta$
 is  bounded by $\sum_{k\in\integer^n}\d^2 |k|^{-n-3}\leq c_n \d^2$.

\nl
Next, recall that properties (P2) and (P3) concern only a {\sl finite} number of $k$,
i.e.,  $k\in\integer^n_*,\  |k|\leq  K_s(\d)$.

\nl 
To show that the set of potentials that do not satisfy (P2) has $\m_s$-measure zero it is enough to check that, for every $k\in\integer^n_*,\  |k|\leq  K_s(\d)$, the set 
${\cal E}^{(k)}$ 
of $f$'s
for which\footnote{Recall the definition of $F_k$ in  \equ{dec}.} $F_k$ has a degenerate critical point has zero $\m_s$-measure. \\
Fix $k\in\integer^n_*,\  |k|\leq  K_s(\d)$ and 
denote points in  $\cE^{(k)}$ by $(\zeta,\f)$, where $\zeta=f_k$ and $\f=\{f_h\}_{h\neq k}$.
Write
$$
F_k(\x)=
\zeta e^{\ii \x} + \bar\zeta e^{-\ii \x} + G(\x)\,,\quad
{\rm where}\ \ 
\zeta:= f_k\ \ {\rm and}\ \ 
 G(\x):=\sum_{|j|\geq 2} f_{jk}e^{\ii j\x}\,.
$$
Now, one checks immediately that $F_k'(\x_0)=0=F_k''(\x_0)$ is equivalent to 
$\zeta=\zeta(\x_0;\f)=\frac12 e^{-\ii \x_0} \big( \ii G'(\x_0)+G''(\x_0) \big)$, which, as $\x_0$ varies in $\torus$, describes a smooth closed ``critical'' curve in $\complex$;
as a side remark, notice that  $\zeta$ depends on $\f$ only through the Fourier coefficients $f_{jk}$ with $|j|\ge 2$.  Thus the section $\cE^{(k)}_\f=\{\zeta\in D: (\zeta,\f)\in\cE^{(k)}\} $ is (a piece of) a smooth curve in $D=\{z\in\complex: |z|\le 1\}$; hence meas$(\cE^{(k)}_\f)=0$ for every $\f$ and by Fubini's theorem $\m_s (\cE^{(k)})=0$, as claimed.

\nl
An analogous  result\footnote{In this case  the critical curve is given by
$\{
\zeta=(-b(\x)\pm \sqrt{b^2(\x)-c(\x)}+\ii G'(\x)) e^{-\ii \x}/2, \ \x\in\real\,, b^2(\x)\geq c(\x)\},$
where $b(\x):=(G''''(\x)-G''(\x))/2$ and $c(\x):=-G''(\x) G''''(\x) +5(G'(\x)+G'''(\x))^2/3.$} 
holds true  for 
(P3). 

\nl
We now show that $\cP_s$ contains an open subset $\cP_s'$ which is dense in the unit ball of  $ \cB_s^n$.

\nl
Let us define $\cP_s'$  as $\cP_s$ but with the difference that  (P1) is replaced by the
{\sl stronger condition}\footnote{Note that $\m_s(\cP'_s)=0$.}

\nl
(P$1'$)
$\exists\, \d>0$ s.t. 
$\displaystyle{  |f_k|\geq\d\ e^{-|k|s} \,,\ \ 
\forall\,k\in\integer^n_*,\  |k|> K_s(\d)}$.

\nl
Let us first prove that $\cP_s'$ is open. Let $f\in\cP_s'.$
We have to show that there exists $\rho>0$ such that if $|g|_s< \rho,$ then $f+g\in\cP_s'.$
Fix $\d>0$ such that (P$1'$) holds and choose $\rho<\d$ small enough such that
$
[K_s(\d)] > K_s(\d')-1\,,$ where $ \d':=\d-\rho
$
and $[\cdot]$ denotes integer part.
Then, it is immediate to verify that $|k|> K_s(\d) \iff |k|> K_s(\d').$
Moreover
$$
|f_k+g_k|e^{|k|s}\geq |f_k| e^{|k|s} -|g|_s \geq \d-\rho=\d'\,, \qquad
\forall\,k\in\integer^n_*,\  |k|> K_s(\d')\,,
$$
namely $f+g$ satisfies  (P$1'$) (with $\d'$ instead of $\d$).
Since (P2) and (P3) are ``open'' conditions and regard only a finite  number of $k$
it is simple to see that they are satisfied also by $f+g$ for $\rho$ small enough.
Then $f+g\in\cP_s'$  for $\rho$ small enough.

\nl
Let us now show that $\cP_s'$  is dense in the unit ball of  $ \cB_s^n$.
Take $f$ in the unit ball of  $ \cB_s^n$ and $0<\theta<1.$ We have to find $\tilde f\in\cP_s'$
with $|\tilde f-f|_s\leq \theta.$
Let $\d:=\theta/4$ and denote by  $f_k$ and $\tilde f_k$ (to be defined) be the Fourier coefficients of, respectively,  
$f$ and $\tilde f$. We, then, let $\tilde f_k=f_k$ unless one of the following two cases occurs: 

\begin{itemize}
\item $k\in\integer^n_*$, $|k|> K_s(\d)$ and $|f_k|e^{|k|s}< \d$, in which case, $\tilde f_k= \delta e^{-|k|s}$;

\item $k\in\integer^n_*$, $|k|\leq K_s(\d)$ and  $F_k$ (defined as in \eqref{dec}) does not satisfy either (P2) or (P3), in which
case, $\tilde f_k$ is chosen at a distance less than $\theta e^{-|k|s}$ from $f_k$  but outside the critical curves defined above.

\end{itemize} 
At this point, it is easy to check that  $\tilde f\in \cP'_s$ and is $\theta$--close to $f$. 

\nl
We finally prove that $\cP_s$ is prevalent.
Consider the following compact subset of
$\ell_\io^n$: 
let $\mathcal K:=\{ z=\{z_k\}_{k\in \integer^n_\sharp  } : z_k\in D_{1/|k|} \},$
where $D_{1/|k|}:=\{w\in\complex:\ |w|\le 1/|k|\},$ 
and let  $\nu$ be 
the unique probability measure supported on $\mathcal K$ such that,
given  Lebesgue measurable sets 
$E_k\subseteq D_{1/|k|}$,
with $E_k\neq D_{1/|k|}$
only for finitely many $k$, one has
$$
\nu \Big(\prod_{k\in\Z^n_\sharp} E_k\Big):=
\prod_{\{k\in \integer^n_\sharp:\, E_k\neq D_{1/|k|}\}} 
\frac{|k|^2}{\pi}{\rm meas}(E_k)\,.
$$ 
The isometry $j_s$ in \eqref{miserere}
naturally induces a  probability measure
$\nu_s$ on  $\cB^n_s$ with support in the compact set $\mathcal K_s
:=j_s^{-1}\mathcal K$.
Now, for $\d>0,$ let $\cP_{s,\d}$ be the set of $f$'s
 in the unit ball of $ \cB_s^n$
satisfying (P1)--(P3), so that $\cP_s=\cup_{\d>0}
\cP_{s,\d}.$
Reasoning as in the proof of 
$\mu_s(\cP_s)=1,$ 
one can show that
$\nu_s(\cP_{s,\d})
\geq 1-{\rm const}\, \d^2.$
It is also easy to check that, for every $g\in  \cB_s^n$, the translated 
set $\cP_{s,\d}+g$ satisfies
$\nu_s(\cP_{s,\d}+g)\geq \nu_s(\cP_{s,\d})$.
Thus,  one gets 
$\nu_s(\cP_s+g)= \nu_s(\cP_s)
= 1$, $\forall\, g\in  \cB_s^n$, which means  that
$\cP_s$ is prevalent
(recall footnote \ref{nurzia}).
 
\Giu
\small
{\bf Acknowledgment.}
We are indebted with V. Kaloshin and A. Sorrentino for useful discussions.

\end{document}